\documentclass[final,5p,times]{elsarticle}

\usepackage{amsmath,amssymb}
\usepackage{booktabs}
\usepackage{graphicx}
\usepackage{array}
\usepackage{placeins}
\usepackage{microtype}
\usepackage{url}
\usepackage{hyperref}
\hypersetup{hidelinks}
\journal{}
\biboptions{numbers,sort&compress}

\begin{document}

\begin{frontmatter}

\title{Boundary-Aware Quantization: Finite-Scale Decision Geometry of Neural Classifiers}

\author[inst1]{O.M. Kiselev\corref{cor1}}
\ead{o.kiselev@innopolis.ru}
\cortext[cor1]{Corresponding author.}

\affiliation[inst1]{organization={Innopolis University},
            city={Innopolis},
            state={Republic of Tatarstan},
            country={Russian Federation}}

\begin{abstract}
We measured quantization-induced decision-boundary changes using local logit-margin radii, first-order boundary displacement, normal variation, slice-boundary Jaccard distance, grid prediction changes, multiclass junction counts, and low-margin boundary-band flips. On the digits benchmark, 8-bit weight quantization preserved all test labels while producing boundary-mask Jaccard \(0.428\) on the PCA slice; at 4 bits, accuracy remained \(0.9733\), while boundary Jaccard rose to \(0.970\) and median local boundary shift reached \(0.0290\). Interpolation between adjacent quantization levels localized the visible reconfigurations at multiclass junctions, with 12, 34, and 17 triple-junction cells in the selected transitions. Calibration-to-test stopping reduced the digits held-out flip rate from \(0.0094\) to \(0.0022\) and boundary Jaccard from \(0.825\) to \(0.524\); the same stopping rule also reduced flips on MNIST and Fashion-MNIST. On official CIFAR-10 subsets, PTQ-W selected by accuracy gave 6-bit flip \(0.0367\) and boundary Jaccard \(0.184\), whereas boundary-aware stopping selected 8-bit flip \(0.0083\) and boundary Jaccard \(0.048\). On full CIFAR-10 with three seeds, 6-bit PTQ-W lost \(0.0029\) accuracy relative to float, changed \(5.3\%\) of held-out decisions, and changed \(24.5\%\) of low-margin boundary-band decisions. A fixed-bit boundary-gap rounding term changed the trade-off at 4 bits by reducing boundary Jaccard from \(0.457\) to \(0.435\) and boundary-band pair-order flip from \(0.3600\) to \(0.3558\), with an accuracy trade-off; the 3-bit stress test exposed the tuning limit of this surrogate. Calibration boundary Jaccard predicted held-out boundary Jaccard across PTQ-W and optimized rounding variants with \(r=0.947\)--\(0.994\).
\end{abstract}


\end{frontmatter}

\section{Introduction}\label{introduction}

Quantization is a standard tool for deploying neural classifiers under
memory, energy, and latency constraints. A full-precision parameter
vector \(\theta\) is replaced by a representable vector
\(\theta_Q=Q(\theta)\), and the resulting model is judged by compression
ratio, runtime, and test accuracy. This evaluation is natural when final
accuracy is the primary target.

In many pattern-recognition pipelines, however, the full-precision model
is already a reference system. The compressed model is then expected to
remain accurate and to behave similarly near class boundaries. Two
classifiers may have almost identical accuracy while assigning different
labels to near-boundary samples, changing local decision regions, or
reorganizing multiclass junctions. Such changes matter when label
consistency across model updates is required, when downstream decisions
depend on stable local predictions, or when adversarial and
near-boundary behavior is part of the acceptance criterion.

Let \(z(x;\theta)\in\mathbb R^K\) be the logit vector and \[
f_\theta(x)=\arg\max_k z_k(x;\theta).
\] Here \(x\) denotes an input, \(K\) is the number of classes,
\(z_k(x;\theta)\) is the logit of class \(k\), and \(f_\theta(x)\) is
the predicted class index for parameters \(\theta\). The class regions
are separated by pairwise boundaries \[
\mathcal B_{ab}(\theta)=\{x:z_a(x;\theta)=z_b(x;\theta)\}.
\] Here \(a\) and \(b\) are distinct class indices, and
\(\mathcal B_{ab}\) is the set where their logits are tied. Quantization
perturbs the logits through \(\theta_Q-\theta\). Away from the boundary,
a large margin may absorb this perturbation. Near a boundary, the same
perturbation can move the separating set across a sample or change how
class regions are connected. Boundary diagnostics reveal this geometry.

This paper asks a practical question: can boundary behavior be measured
on a calibration set and then used to decide when quantization should
stop? We answer this by proposing a finite-scale boundary audit for
quantized classifiers. On a fixed two-dimensional slice and grid, we
measure observed boundary masks, grid prediction changes, and multiclass
junction cells. These quantities are deployment diagnostics for
reference consistency at the monitored finite scale.

The work is positioned at the intersection of quantization, robustness,
and decision-boundary analysis. Neural network quantization has been
developed for efficient low-precision inference, including integer
deployment and quantization-aware training
\citep{jacob2018quantization, krishnamoorthi2018quantizing}. Data-aware
post-training methods, such as AdaRound, improve rounding decisions by
approximating local loss changes \citep{nagel2020up}. These methods
primarily optimize task loss, accuracy, or deployment cost.

The adversarial robustness literature shows that high accuracy can
coexist with geometric instability. Small input perturbations can cross
a decision boundary and change predictions
\citep{szegedy2014intriguing, goodfellow2015explaining}.
Boundary-distance and decision-geometry analyses provide the closest
methodological background: DeepFool uses a local linear approximation to
estimate boundary-crossing perturbations, curvature-based work relates
robustness to boundary curvature, and empirical decision-boundary
studies inspect boundary shapes in input space
\citep{moosavi2016deepfool, moosavi2019curvature, mickisch2020understanding}.
Quantization and robustness have also been studied together: Defensive
Quantization controls Lipschitz behavior during quantization
\citep{lin2019defensive}, and other work evaluates or optimizes
adversarial robustness of quantized networks
\citep{gorsline2021adversarial, song2021layerwise}.

The boundary-band analysis used below is a reference-preservation
diagnostic. The low-margin band localizes regressions introduced by
compression relative to a validated full-precision model. In this
setting, the relevant events include errors against the ground-truth
label, changes of the reference decision, and changes of the reference
top-2 class ordering. Thus boundary-band flip is a stratified regression
metric tied to quantization-induced decision-geometry drift.

The rest of the paper follows this structure. Section 2 combines local
boundary perturbation estimates with finite-scale slice metrics,
multiclass-node retention, boundary-aware stopping, and the
boundary-aware quantization objective. Section 3 reports the experiments
on digits, MNIST, Fashion-MNIST, and official CIFAR-10 with
post-training quantization, activation quantization, optimized rounding,
and fake-quant quantization-aware training; it also records
reproducibility and data availability. Section 4 gives the conclusions.

The main empirical message is that quantization can preserve or improve
accuracy while changing the reference decision geometry. Boundary-aware
stopping serves as an audit layer that tests whether a chosen
quantization method preserves the behavior of the reference classifier
at the finite scales being monitored.

\section{Boundary Change and Finite-Scale
Metrics}\label{boundary-change-and-finite-scale-metrics}

For two classes \(a,b\), define the logit gap \[
g_{ab}(x;\theta)=z_a(x;\theta)-z_b(x;\theta).
\] The corresponding boundary is the zero set \(g_{ab}=0\). In this
notation, \(a\) and \(b\) are the two compared class indices, and a
positive value of \(g_{ab}\) means that class \(a\) has a larger logit
than class \(b\). For a labelled sample \((x,y)\), define the true-label
margin \[
m_y(x;\theta)=z_y(x;\theta)-\max_{r\ne y}z_r(x;\theta).
\] Here \(y\) is the ground-truth label and \(r\) ranges over all
competing labels. If the uniform logit perturbation satisfies
\(\|z(x;\theta_Q)-z(x;\theta)\|_\infty\le \eta\) and
\(m_y(x;\theta)>2\eta\), then the quantized classifier keeps the correct
label at \(x\). Thus samples with large logit margin are stable, while
near-boundary samples are potentially sensitive to quantization.

When \(g_{ab}\) is differentiable and \(\nabla_x g_{ab}(x;\theta)\ne0\),
the first-order distance from \(x\) to the pairwise boundary is \[
\rho_{ab}(x;\theta)=
\frac{|g_{ab}(x;\theta)|}{\|\nabla_x g_{ab}(x;\theta)\|_2}.
\] Here \(\nabla_x\) denotes the gradient with respect to the input and
\(\|\cdot\|_2\) is the Euclidean norm. The quantity \(\rho_{ab}\) is the
first-order normal distance from \(x\) to the pairwise boundary
\(\mathcal B_{ab}\). The quantization-induced displacement of this local
boundary is estimated by \[
d^Q_{ab}(x)\approx
-
\frac{g_{ab}(x;\theta_Q)-g_{ab}(x;\theta)}
{\|\nabla_x g_{ab}(x;\theta)\|_2}.
\] The scalar \(d^Q_{ab}(x)\) is a signed first-order displacement along
the reference normal. The ratio \(|d^Q_{ab}(x)|/\rho_{ab}(x;\theta)\) is
a direct warning signal: if it is large, the quantized boundary can
cross the sample or change the local decision geometry around it. We
also measure the cosine similarity between \(\nabla_x g_{ab}(x;\theta)\)
and \(\nabla_x g_{ab}(x;\theta_Q)\), since a boundary may rotate while
the nearby prediction remains unchanged.

The same mechanism follows from the implicit function theorem. If \(g\)
and \(g_Q\) are \(C^1\), \(\|\nabla g\|\ge\mu>0\) near a compact regular
boundary segment, and \(\|g_Q-g\|_{C^0}\le\varepsilon\), then the
quantized zero set remains locally within distance
\(O(\varepsilon/\mu)\). Here \(\mu\) is a lower bound on the reference
boundary gradient and \(\varepsilon\) is a uniform functional
perturbation bound. Low-gradient boundary regions are therefore more
sensitive to quantization-induced functional perturbations. These local
estimates are combined with finite-scale slice diagnostics.

The full decision boundary is high-dimensional, so we observe it on a
prescribed family of two-dimensional affine slices \(\Pi\). On each
slice, a regular grid is labelled by the classifier. A grid cell is a
boundary cell if its vertices include more than one predicted label. The
resulting boundary mask is denoted \(S_\Pi(\theta)\).

The primary slice metric is the boundary-mask Jaccard distance \[
J_\Delta(S,S_Q)=1-\frac{|S\cap S_Q|}{|S\cup S_Q|}.
\] Here \(S=S_\Pi(\theta)\) is the full-precision boundary mask,
\(S_Q=S_\Pi(\theta_Q)\) is the quantized boundary mask, \(|\cdot|\)
denotes the number of grid cells in a set, and \(\Delta\) indicates that
the comparison is made on the fixed slice grid. The metric measures
boundary reconfiguration relative to the full-precision reference. We
also report grid prediction flip rate on the same slices, because a
boundary mask can move substantially while affecting part of the
monitored slice. These slice quantities are reference-consistency
diagnostics at the chosen resolution.

\subsection{Multiclass Junctions}\label{multiclass-junctions}

Multiclass boundaries contain junctions where three or more class
regions meet. Such nodes are more sensitive than ordinary pairwise arcs
because several logit equalities must hold simultaneously. A standard
local model is the max-affine boundary of the upper envelope of affine
logits \citep{richter2005first, maclagan2015introduction}. In the plane,
stable vertices are typically trivalent, while four-valent and
higher-valent nodes are degenerate.

The simplest four-ray model is \[
M_\delta(u,v)=\max\{u,\;v,\;-u,\;-v+\delta\}.
\] Here \(u\) and \(v\) are local coordinates on the slice, \(M_\delta\)
is the upper envelope of four affine logits, and \(\delta\) is a scalar
perturbation of one branch. At \(\delta=0\), the four regions meet at
the origin. For \(\delta>0\), this node splits into two trivalent nodes
\[
p_-=\left(-\frac{\delta}{2},\frac{\delta}{2}\right),
\qquad
p_+=\left(\frac{\delta}{2},\frac{\delta}{2}\right),
\] joined by the short edge \(v=\delta/2\), \(|u|\le\delta/2\). This
standard model gives the diagnostic mechanism: quantization can break
ties among competing logits and resolve a high-valence node into a small
cluster of lower-valence nodes.

This motivates a finite-scale node-retention criterion. Let \(x_0\) be a
reference junction on a slice, and let \(I=\{i_1,\ldots,i_k\}\) be the
participating classes. Define \[
G_I(x;\theta)=
\bigl(
z_{i_2}(x;\theta)-z_{i_1}(x;\theta),
\ldots,
z_{i_k}(x;\theta)-z_{i_1}(x;\theta)
\bigr).
\] The vector \(G_I\) contains the \(k-1\) logit tie equations that
define a junction with class set \(I\), using class \(i_1\) as the
reference logit. For tolerance \(\tau\), the node-retention radius is \[
R_I^\tau(x_0;\theta_Q)=
\inf\{r\ge0:\exists x\in B_\Pi(x_0,r)
\text{ with }\|G_I(x;\theta_Q)\|_\infty\le\tau\}.
\] Here \(B_\Pi(x_0,r)\) is the radius-\(r\) ball inside the slice
\(\Pi\), \(\|\cdot\|_\infty\) is the maximum absolute component, and
\(\tau\) is the allowed logit-tie tolerance. In a discrete
implementation, the same idea is measured by the enclosing radius of the
node cluster or by the anchored radius of the cells with the same
class-incidence pattern. The practical requirement is that the resolved
cluster remains inside an admissible ball.

\subsection{Boundary-Aware
Stopping}\label{boundary-aware-stopping}

Given bit-widths \(b_1>b_2>\cdots>b_m\), boundary-aware stopping
evaluates each quantized model on a calibration set and a prescribed
slice family. It selects the smallest bit-width satisfying task-specific
constraints: \[
\Delta\operatorname{Acc}(b_j)\le\tau_{\mathrm{acc}},
\qquad
\operatorname{Flip}_{cal}(b_j)\le\tau_{\mathrm{flip}},
\] \[
\operatorname{median}_x
\frac{|d^{Q_{b_j}}(x)|}{\operatorname{median}_x\rho(x;\theta)}
\le\tau_d,
\qquad
\operatorname{median}_x \cos(n_\theta(x),n_{\theta_Q}(x))\ge\tau_n,
\] and, when boundary preservation is required, \[
J_\Delta(S_\Pi(\theta),S_\Pi(\theta_{Q_{b_j}}))\le\tau_J
\] on the monitored slices. Here \(\theta_{Q_{b_j}}\) denotes the model
quantized to \(b_j\) bits, \(\Delta\operatorname{Acc}\) is the
validation accuracy loss relative to the full-precision reference,
\(\operatorname{Flip}_{cal}\) is the calibration label-flip rate
relative to the reference, \(d^{Q_{b_j}}\) is the corresponding top-2
boundary displacement, \(\rho\) is the corresponding top-2 local radius,
\(n_\theta(x)=\nabla_x g(x;\theta)/\|\nabla_x g(x;\theta)\|_2\) is the
local unit normal of the reference top-2 boundary, \(g\) is the top-2
logit gap at \(x\), and
\(\tau_{\mathrm{acc}},\tau_{\mathrm{flip}},\tau_d,\tau_n,\tau_J\) are
deployment tolerances. Optional node constraints can bound node
displacement, node-cluster radius, or changes in the number of junction
cells. If all tested quantized candidates violate the constraints, the
protocol reports a constraint violation together with the least
violating candidate, giving the user a choice between relaxing the
tolerances, changing the quantizer, or accepting the measured boundary
change.

\subsection{Boundary-Aware Quantization
Objective}\label{boundary-aware-quantization-objective}

The same quantities can also be used inside the quantization procedure.
Let \(\mathcal G_b\) be the set of parameters representable at bit-width
\(b\), and let \(\mathcal C=\{c_1,\ldots,c_M\}\) be a finite family of
monitored boundary objects selected on calibration slices. These objects
may include ordinary pairwise boundary cells, multiclass junctions, or
high-valence node clusters. Each object has an admissible radius
\(r_c>0\) and weight \(w_c\ge0\). A boundary-aware quantizer can then be
written as the constrained problem \[
\begin{aligned}
\min_{\theta_Q\in\mathcal G_b}\quad
&\mathcal L_{\mathrm{cal}}(\theta_Q)
+\lambda_\theta\|\theta_Q-\theta\|^2
+\lambda_BD_{\mathrm{bdry}}(\theta_Q,\theta)\\
\text{subject to}\quad
&R_c(\theta_Q)\le r_c,\quad c=1,\ldots,M .
\end{aligned}
\] Here \(\mathcal L_{\mathrm{cal}}\) is the calibration loss,
\(\lambda_\theta\) and \(\lambda_B\) are nonnegative weights,
\(D_{\mathrm{bdry}}\) is a boundary discrepancy functional, and \(R_c\)
is the measured displacement or cluster radius of object \(c\). The
coefficient \(w_c\) is a user-chosen priority weight for the monitored
object \(c\). It can encode the importance of the classes meeting at
this boundary object, the density of calibration samples near it, the
reliability of its detection on the slice grid, or the application cost
of changing decisions near that object. In the unweighted case \(w_c=1\)
for all monitored objects; in a normalized implementation one may use
\(\sum_{c\in\mathcal C}w_c=1\). Larger \(w_c\) allocates more of the
quantization objective to preserving object \(c\) relative to other
monitored arcs and nodes.

For degenerate nodes, the object \(c\) has a reference location \(x_c\),
a participating label set \(I_c\), an audit search radius \(\rho_c\),
and a logit-tie tolerance \(\tau_c\). The split cluster generated by
quantization is \[
\mathcal N_c(\theta_Q)=
\{x\in B_\Pi(x_c,\rho_c):\|G_{I_c}(x;\theta_Q)\|_\infty\le\tau_c\}.
\] Two useful sizes of this cluster are the anchored radius \[
R_c^{\mathrm{anc}}(\theta_Q)=
\sup_{x\in\mathcal N_c(\theta_Q)}\|x-x_c\|_2
\] and the minimum enclosing radius \[
R_c^{\mathrm{enc}}(\theta_Q)=
\inf_{y\in\Pi}\sup_{x\in\mathcal N_c(\theta_Q)}\|x-y\|_2 .
\] The anchored radius measures how far the resolved nodes move from the
original degenerate node, while the enclosing radius measures the
smallest ball on the slice that contains the whole resolved cluster. In
a grid implementation, \(\mathcal N_c\) is the finite set of cells or
grid points satisfying the same incidence/tie condition, and the suprema
and infimum are computed over that finite set. The empty-cluster
indicator
\(E_c(\theta_Q)=\mathbf 1[\mathcal N_c(\theta_Q)=\varnothing]\) records
disappearance of the monitored incidence pattern.

This gives a direct optimization problem for keeping resolved nodes
inside small balls: \[
\begin{aligned}
\min_{\theta_Q\in\mathcal G_b}\quad
&\mathcal L_{\mathrm{cal}}(\theta_Q)
+\lambda_\theta\|\theta_Q-\theta\|^2
+\lambda_A\sum_{c\in\mathcal C_N}w_c
\left(\frac{R_c^{\mathrm{anc}}(\theta_Q)}{r_c}\right)^2
\\
&+\lambda_E\sum_{c\in\mathcal C_N}w_c
\left(\frac{R_c^{\mathrm{enc}}(\theta_Q)}{r_c}\right)^2
+\lambda_0\sum_{c\in\mathcal C_N}w_cE_c(\theta_Q),
\end{aligned}
\] where \(\mathcal C_N\subseteq\mathcal C\) is the subset of monitored
degenerate nodes and \(\lambda_A,\lambda_E,\lambda_0\ge0\) control
anchored displacement, cluster diameter, and disappearance penalties. A
minimax version emphasizes the worst preserved node: \[
\begin{aligned}
\min_{\theta_Q\in\mathcal G_b}\quad
&\mathcal L_{\mathrm{cal}}(\theta_Q)
+\lambda_\theta\|\theta_Q-\theta\|^2
+\lambda_M\max_{c\in\mathcal C_N} w_c\\
&\times
\frac{\max\{R_c^{\mathrm{anc}}(\theta_Q),
R_c^{\mathrm{enc}}(\theta_Q)\}}{r_c}
+\lambda_0\sum_{c\in\mathcal C_N}w_cE_c(\theta_Q).
\end{aligned}
\] This formulation covers general \(k\)-valent nodes: the vector
\(G_{I_c}\) contains \(k-1\) tie equations, and the cluster
\(\mathcal N_c\) represents the resolved finite-grid trace of that
degenerate incidence pattern. It also makes the engineering problem
explicit: quantization must distribute rounding error across many
monitored arcs and nodes, so a locally harmless rounding choice can
still enlarge another node cluster.

Equivalently, a soft version minimizes \[
\begin{aligned}
\min_{\theta_Q\in\mathcal G_b}\quad
&\mathcal L_{\mathrm{cal}}(\theta_Q)
+\lambda_\theta\|\theta_Q-\theta\|^2
+\lambda_BD_{\mathrm{bdry}}(\theta_Q,\theta)\\
&+\lambda_N\sum_{c\in\mathcal C}w_c
\left[\frac{R_c(\theta_Q)}{r_c}-1\right]_+^2 .
\end{aligned}
\] over the quantization grid. Here \(D_{\mathrm{bdry}}\) may combine
boundary-mask Jaccard distance, grid prediction changes, normal
variation, and node-retention penalties on the monitored slice family.
This formulation makes explicit that preserving many boundary objects is
a global allocation problem: a rounding choice that is harmless for one
node may move another node or boundary arc.

The experiments below use a differentiable surrogate of this objective
for optimized rounding: near-boundary calibration points are selected by
small full-precision top-2 margin, and the rounding loss penalizes
changes in their top-2 logit gaps. Exact boundary-mask distances and
node-cluster radii are evaluated as finite-grid audit metrics. This
ablation tests whether boundary information is usable by a quantizer and
clarifies the accuracy-boundary trade-off.

\section{Experiments}\label{experiments}

We evaluate whether boundary metrics reveal quantization effects that
are hidden by accuracy and whether calibration metrics predict held-out
behavior. The implementation produces raw CSV/JSON files and figures
under \texttt{results/}. The experiments include a controlled MLP on
\texttt{sklearn\ digits}, stratified MNIST and Fashion-MNIST subsamples,
official CIFAR-10 subsets with a small CNN, a full-dataset reduced
residual CNN check, post-training weight quantization (PTQ-W),
post-training weight-and-activation quantization (PTQ-WA), and
fake-quant quantization-aware training with logit distillation (QAT-WA).

For a tensor \(T\), the uniform affine fake-quantizer used in the
experiments has the form \[
\widehat T=s\left(
\operatorname{clip}\left(\operatorname{round}(T/s)+z_0,q_{\min},q_{\max}\right)-z_0
\right).
\] Here \(b\) is the bit-width, \(q_{\min}\) and \(q_{\max}\) are the
integer endpoints of the \(b\)-bit quantization range, \(s>0\) is the
scale, \(z_0\) is the integer zero point, and \(\widehat T\) is the
dequantized tensor used by the floating-point simulator. In the
symmetric weight case, \(z_0=0\) and \(s\) is determined from the tensor
range.

PTQ-W means post-training weight quantization. The full-precision model
is trained first, its convolutional and linear weights are mapped to the
\(b\)-bit grid, and the quantized weights are evaluated with activations
computed in floating point. Biases and batch-normalization parameters
are kept in full precision in these experiments.

PTQ-WA means post-training weight-and-activation quantization. The
weights are quantized as in PTQ-W, and activation tensors produced
during the forward pass are fake-quantized to the same bit-width. This
family measures the additional boundary drift caused by activation
quantization after training.

QAT-WA means quantization-aware training with fake quantization of
weights and activations. Starting from the full-precision model,
fake-quantization operators are inserted into the forward pass and the
model is fine-tuned under the quantized computation graph. The loss
combines supervised cross-entropy with logit distillation from the
full-precision reference, so QAT-WA optimizes task performance while the
low-precision simulation is active.

All slice metrics are finite-grid quantities. PCA slices expose
high-variance data directions, while random slices reduce dependence on
one coordinate system.

The computational overhead is controlled by the number of monitored
samples, slices, and grid points. The protocol is intended as an offline
calibration audit for candidate quantized models. Point margins, flips,
and boundary-band flips require forward passes on calibration or
held-out samples; local displacement and normals additionally require
input gradients for the reference top-2 gaps. Slice masks require
\(Sg^2\) grid evaluations for \(S\) monitored slices of grid width
\(g\). In the CIFAR-10 runs we use \(S=3\), \(g=16\), 2,000 monitored
samples per split for point metrics, and 4,096 calibration images for
optimized rounding.

\subsection{Controlled Digits
Experiment}\label{controlled-digits-experiment}

The digits model is a one-hidden-layer tanh network. We quantize weights
and keep biases in full precision. Uniform affine per-tensor
quantization is tested at 8, 4, 3, and 2 bits.

\begin{figure*}[p]
\centering
\begin{minipage}[t]{0.50\textwidth}
\centering
\scriptsize
\setlength{\tabcolsep}{1.6pt}
\renewcommand{\arraystretch}{1.06}
\begin{tabular}{lrrrrrrr}
\hline
model & acc. & flip & q10 \(m_y\) & q10 \(\rho\) & med. \(d\) & q90 \(d\) & med. cos \\
\hline
float & 0.9756 & 0.0000 & 2.8433 & 0.1808 & 0.0000 & 0.0000 & 1.0000 \\
8-bit & 0.9756 & 0.0000 & 2.8499 & 0.1797 & 0.0020 & 0.0051 & 1.0000 \\
4-bit & 0.9733 & 0.0022 & 3.0195 & 0.1853 & 0.0290 & 0.0621 & 0.9946 \\
3-bit & 0.9644 & 0.0222 & 1.9290 & 0.1305 & 0.0887 & 0.1992 & 0.9745 \\
2-bit & 0.8778 & 0.1133 & -0.3377 & 0.0450 & 0.1951 & 0.4689 & 0.8421 \\
\hline
\end{tabular}
\end{minipage}\hfill
\begin{minipage}[t]{0.34\textwidth}
\centering
\vspace{-3.1em}
\includegraphics[width=.82\linewidth]{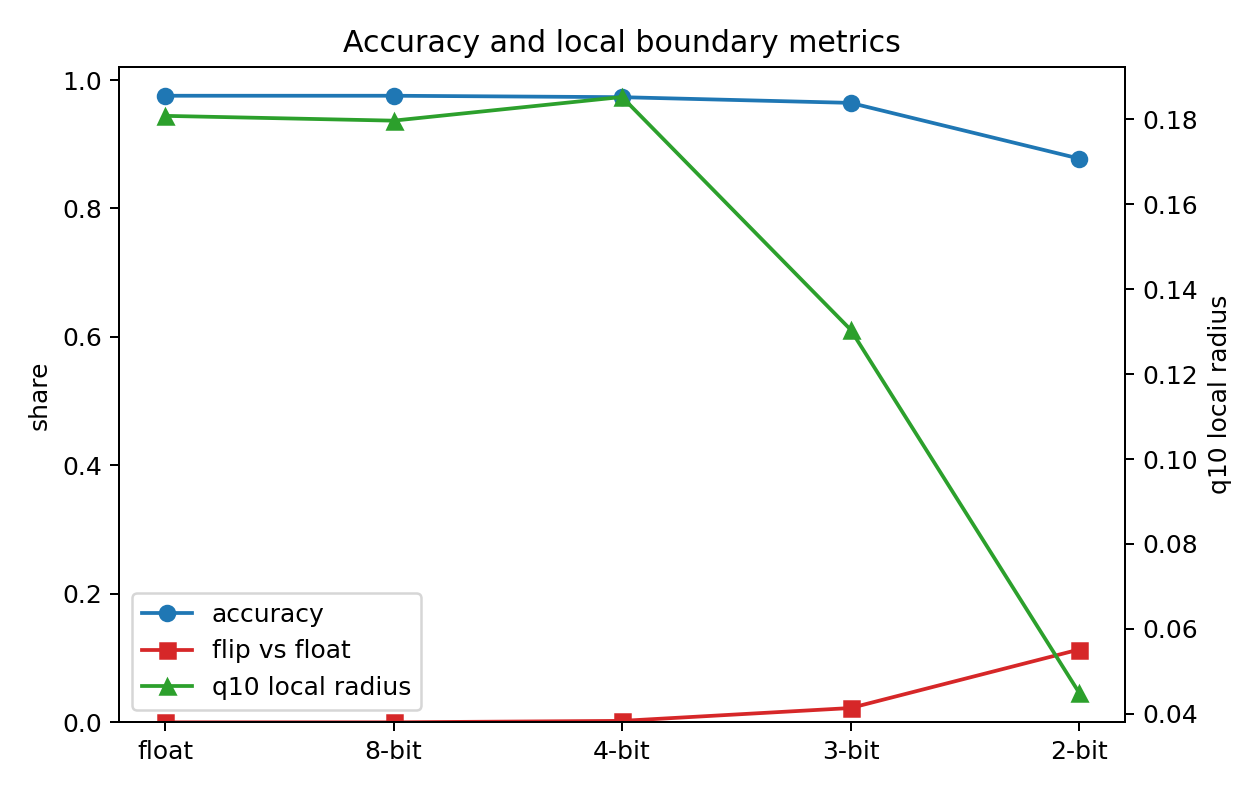}
\end{minipage}
\caption{Controlled digits quantization trade-off.}
\label{fig:tradeoff-prl}
\end{figure*}

The 8-bit model preserves all test labels. At 4 bits, accuracy remains
close to the float model, while the local boundary shift is already
visible. At 3 and 2 bits, the radius quantile drops and boundary
displacement increases strongly. This supports the use of local geometry
as an early warning signal.

On the PCA slice, the boundary-mask Jaccard distance grows sharply even
when accuracy is still high.

\begin{figure*}[p]
\centering
\scalebox{1}[0.88]{\includegraphics[width=.92\textwidth]{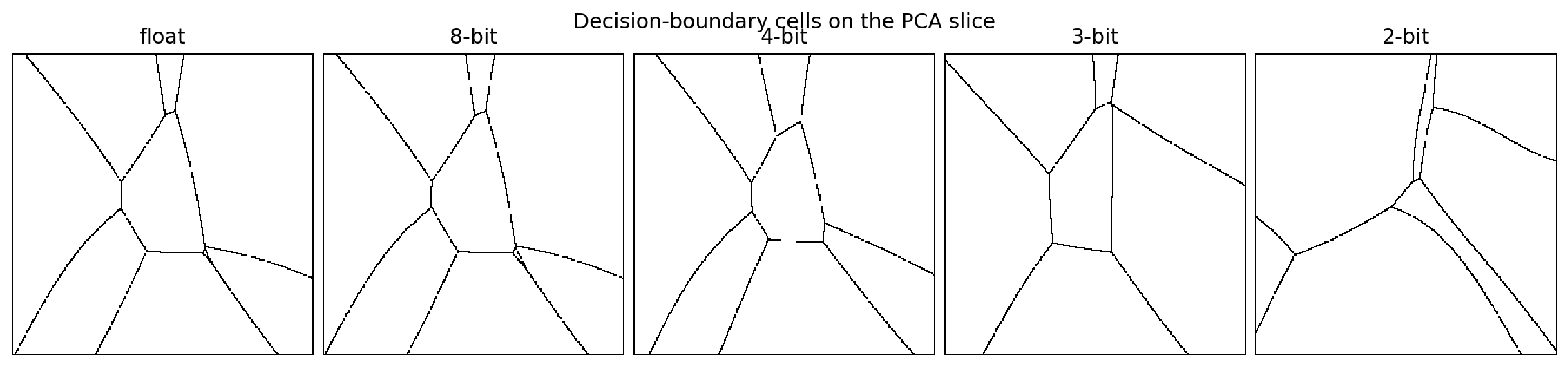}}
\caption{Boundary cells on a PCA slice of the input space.}\label{fig:pca-mask-prl}
\end{figure*}

Table~\ref{tab:digits-slice-boundary} shows why boundary Jaccard and grid flip should be read
together. From float to 4-bit, the boundary mask changes substantially,
while part of the monitored grid changes label. The dominant event is
reference-boundary reconfiguration.

\subsection{Interpolation and Junction
Reconfiguration}\label{interpolation-and-junction-reconfiguration}

To interpret visible boundary changes, we interpolate between
neighbouring quantized parameter vectors: \[
\theta(t)=(1-t)\theta_A+t\theta_B,\qquad t\in[0,1].
\] Here \(\theta_A\) and \(\theta_B\) are neighbouring quantized
parameter vectors and \(t\) is the interpolation parameter. For a
two-class boundary on a slice, a Morse-type event would occur near
\(F_{ab}=0\), \(\partial_uF_{ab}=0\), \(\partial_vF_{ab}=0\). Here
\(F_{ab}(u,v;t)=z_a(x(u,v);\theta(t))-z_b(x(u,v);\theta(t))\) is the
logit gap restricted to slice coordinates \(u,v\) along \(t\). In our
controlled slice, the selected events are localized at multiclass
junction cells.

\begin{table*}[p]
\centering
\begin{minipage}[t]{0.34\textwidth}
\centering
\scriptsize
\setlength{\tabcolsep}{3.2pt}
\renewcommand{\arraystretch}{1.05}
\begin{tabular}{lrr}
\hline
model & grid flip & \(J_\Delta\) \\
\hline
float & 0.0000 & 0.000 \\
8-bit & 0.0056 & 0.428 \\
4-bit & 0.0949 & 0.970 \\
3-bit & 0.3221 & 0.985 \\
2-bit & 0.5933 & 0.986 \\
\hline
\end{tabular}
\caption{Slice-level boundary change on the digits PCA slice.}
\label{tab:digits-slice-boundary}
\end{minipage}\hfill
\begin{minipage}[t]{0.62\textwidth}
\centering
\scriptsize
\setlength{\tabcolsep}{3.2pt}
\renewcommand{\arraystretch}{1.05}
\begin{tabular}{lccccc}
\hline
transition & selected \(t^\ast\) & event class & two-class saddle cells & triple-junction cells & Jaccard step \\
\hline
8-bit \(\to\) 4-bit & 0.288 & multiclass junction & 0 & 12 & 0.122 \\
4-bit \(\to\) 3-bit & 0.988 & multiclass junction & 0 & 34 & 0.405 \\
3-bit \(\to\) 2-bit & 0.700 & multiclass junction & 0 & 17 & 0.487 \\
\hline
\end{tabular}
\caption{Selected interpolation events between adjacent quantization levels.}
\label{tab:interpolation-events}
\end{minipage}
\end{table*}

The detector returns zero two-class diagonal saddle cells at the
selected events. Representative event clusters involve three nearby
labels. Thus the appropriate finite-grid description is multiclass
junction reattachment or reconnection of the boundary graph.

\subsection{Calibration-to-Test
Stopping}\label{calibration-to-test-stopping}

We next evaluate stopping profiles. An accuracy profile accepts the
smallest bit-width whose validation accuracy and flip rate remain within
tolerance. A boundary-aware profile additionally constrains local
displacement, normal stability, and boundary Jaccard distance. The
thresholds represent deployment tolerances.

For the digits benchmark with five random seeds, bit-widths
\(8,6,5,4,3,2\), and both PCA and random slices, accuracy stopping
selects 4-bit quantization, whereas boundary-aware stopping selects
8-bit. This reduces the mean held-out flip rate from \(0.0094\) to
\(0.0022\) and the mean boundary Jaccard distance from \(0.825\) to
\(0.524\). Calibration boundary metrics also have predictive signal on
this benchmark: displacement normalized by baseline radius correlates
with held-out label flip rate at \(r=0.94\), normal cosine correlates
negatively with held-out flip rate at \(r=-0.97\), and calibration
boundary Jaccard predicts held-out boundary Jaccard at \(r=0.85\).

The same pattern appears on stratified MNIST and Fashion-MNIST
subsamples (Table \ref{tab:cross-dataset-stopping}).

\begin{table*}[p]
\centering
\scriptsize
\setlength{\tabcolsep}{4pt}
\renewcommand{\arraystretch}{1.08}
\begin{tabular}{llcrr}
\hline
dataset & profile & bit & flip vs float, mean [95\% CI] & \(J_\Delta\), mean [95\% CI] \\
\hline
Digits & accuracy & 4-bit & 0.0094 [0.0039, 0.0150] & 0.825 [0.720, 0.922] \\
Digits & boundary-aware & 8-bit & 0.0022 [0.0000, 0.0056] & 0.524 [0.372, 0.676] \\
MNIST & accuracy & 3-bit & 0.0347 [0.0280, 0.0413] & 0.425 [0.406, 0.443] \\
MNIST & boundary-aware & 4-bit & 0.0120 [0.0057, 0.0183] & 0.303 [0.243, 0.364] \\
Fashion-MNIST & accuracy & 3-bit & 0.0677 [0.0297, 0.1077] & 0.610 [0.386, 0.833] \\
Fashion-MNIST & boundary-aware & 5-bit & 0.0170 [0.0100, 0.0227] & 0.482 [0.287, 0.688] \\
\hline
\end{tabular}
\caption{Calibration-to-test stopping outcomes on digits, MNIST, and Fashion-MNIST.}
\label{tab:cross-dataset-stopping}
\end{table*}

Across datasets, accuracy stopping tends to choose lower bit-widths and
produces larger held-out flip rates. Boundary-aware stopping is more
conservative and reduces decision changes. Calibration displacement
normalized by baseline radius correlates with held-out flip rate at
\(r=0.94\), \(0.98\), and \(0.96\) on Digits, MNIST, and Fashion-MNIST,
respectively. Calibration normal cosine correlates negatively with
held-out flip rate at \(r=-0.97\), \(-0.96\), and \(-0.98\). Calibration
boundary Jaccard predicts held-out boundary Jaccard at \(r=0.85\),
\(0.97\), and \(0.93\).

\subsection{Official CIFAR-10 Experiments and Quantization
Families}\label{official-cifar-10-experiments-and-quantization-families}

We then test whether the same behavior appears beyond shallow MLPs. On
official CIFAR-10 subsets, we use a small CNN with 10,000 stratified
images, three seeds, bit-widths \(8,6,4,3\), and three quantization
families: PTQ-W, PTQ-WA, and QAT-WA with logit distillation from the
full-precision model.

\begin{table*}[p]
\centering
\scriptsize
\setlength{\tabcolsep}{4pt}
\renewcommand{\arraystretch}{1.08}
\caption{Official CIFAR-10 subset stopping outcomes.}
\label{tab:cifar-stopping}
\begin{tabular}{llccc}
\hline
family & profile & mode selected bit & test flip vs float, mean [95\% CI] & boundary Jaccard, mean [95\% CI] \\
\hline
PTQ-W & accuracy & 6-bit & 0.0367 [0.0332, 0.0435] & 0.184 [0.109, 0.329] \\
PTQ-W & boundary-aware & 8-bit & 0.0083 [0.0080, 0.0090] & 0.048 [0.024, 0.077] \\
PTQ-WA & accuracy & 6-bit & 0.0463 [0.0385, 0.0510] & 0.207 [0.133, 0.281] \\
QAT-WA & accuracy & 4-bit & 0.1668 [0.1603, 0.1734] & 0.480 [0.387, 0.647] \\
\hline
\end{tabular}
\end{table*}

For PTQ-WA and QAT-WA, all tested quantized bit-widths violate the
strict boundary-aware profile, so the table reports the quantized
selections. The full-precision reference has zero flip and zero boundary
Jaccard by definition. This result clarifies the role of QAT. QAT-WA can
preserve or improve classification accuracy while producing a
geometrically different classifier relative to the original reference.
Under accuracy stopping it selects 4-bit quantization with held-out flip
rate \(0.1668\) and boundary Jaccard \(0.480\). Boundary metrics capture
a reference-consistency axis separate from the QAT training objective.

We also ran optimized rounding checks on the same official CIFAR-10
subset scale. The first optimized baseline uses soft rounding variables
for each weight tensor and minimizes a logit-distillation loss to the
full-precision model on calibration images, then projects the variables
to hard quantized weights. This gives an AdaRound-style weight-rounding
baseline. The second optimized baseline adds the proposed boundary term.
The experiment is a fixed-bit ablation of the boundary term. To separate
the effect of the quantizer from the effect of bit-width selection, the
main comparison below is made at fixed bit-widths.

For an engineering safety check, we also report a boundary-band flip
rate. The band contains the \(20\%\) of held-out samples with the
smallest full-precision top-2 logit gap. It measures
reference-regression risk: how often quantization changes the decision
precisely where the reference classifier is closest to a boundary. In
the table, Logit denotes logit rounding-W and Boundary denotes boundary
rounding-W.

\begin{table*}[p]
\centering
\begin{minipage}[t]{0.56\textwidth}
\centering
\scriptsize
\setlength{\tabcolsep}{2.2pt}
\renewcommand{\arraystretch}{1.02}
\begin{tabular}{rlrrrrr}
\hline
bit & method & acc. & flip & band & pair & \(J_\Delta\) \\
\hline
6 & PTQ & 0.5565 & 0.0340 & 0.1692 & 0.1667 & 0.223 \\
6 & Logit & 0.5560 & 0.0328 & 0.1642 & 0.1608 & 0.211 \\
6 & Boundary & 0.5555 & 0.0328 & 0.1642 & 0.1608 & 0.225 \\
4 & PTQ & 0.5495 & 0.1182 & 0.3983 & 0.3708 & 0.462 \\
4 & Logit & 0.5487 & 0.1148 & 0.3883 & 0.3600 & 0.457 \\
4 & Boundary & 0.5467 & 0.1180 & 0.3900 & 0.3558 & 0.435 \\
3 & PTQ & 0.5055 & 0.2620 & 0.5350 & 0.4317 & 0.560 \\
3 & Logit & 0.5067 & 0.2603 & 0.5408 & 0.4433 & 0.564 \\
3 & Boundary & 0.5122 & 0.2617 & 0.5492 & 0.4592 & 0.648 \\
\hline
\end{tabular}
\caption{Fixed-bit CIFAR-10 subset rounding comparison.}
\label{tab:fixed-bit-rounding}
\end{minipage}\hfill
\begin{minipage}[t]{0.41\textwidth}
\centering
\scriptsize
\setlength{\tabcolsep}{2.0pt}
\renewcommand{\arraystretch}{1.02}
\begin{tabular}{lrrrrrr}
\hline
family & bit & acc. & flip & band & pair & \(J_\Delta\) \\
\hline
PTQ-W & 8 & 0.6993 & 0.0159 & 0.0797 & 0.0789 & 0.187 \\
PTQ-W & 6 & 0.6971 & 0.0530 & 0.2453 & 0.2297 & 0.446 \\
PTQ-W & 4 & 0.6196 & 0.2554 & 0.5422 & 0.4318 & 0.746 \\
PTQ-W & 3 & 0.3005 & 0.6959 & 0.8164 & 0.4768 & 0.814 \\
PTQ-WA & 8 & 0.6997 & 0.0176 & 0.0878 & 0.0861 & 0.186 \\
PTQ-WA & 6 & 0.6965 & 0.0613 & 0.2717 & 0.2515 & 0.434 \\
PTQ-WA & 4 & 0.5826 & 0.3109 & 0.5871 & 0.4419 & 0.803 \\
PTQ-WA & 3 & 0.2121 & 0.7878 & 0.8643 & 0.4825 & 0.773 \\
\hline
\end{tabular}
\caption{Full CIFAR-10 reduced-residual scale check.}
\label{tab:full-cifar}
\end{minipage}
\end{table*}

All entries are means over three seeds on the held-out split. The
boundary-band columns show why boundary diagnostics are relevant for
safety-sensitive deployment updates: quantization flips are concentrated
near the full-precision boundary. For example, at 6 bits the overall
PTQ-W flip rate is \(0.0340\), while the flip rate inside the low-margin
boundary band is \(0.1692\). Thus global accuracy and even global flip
rate understate the local regression risk near the reference boundary.

The fixed-bit comparison remains a trade-off analysis. At 6 bits, logit
rounding and boundary-aware rounding have the same band flip and pair
flip, while logit rounding has a lower slice-boundary Jaccard. At 4
bits, boundary-aware rounding gives the lowest boundary Jaccard and the
lowest pair-order flip in the boundary band, with lower accuracy than
PTQ-W and a higher overall flip rate than logit rounding. The 3-bit
regime is an aggressive stress test: boundary-aware rounding improves
accuracy relative to PTQ-W and logit rounding, while boundary Jaccard
and boundary-band pair flip deteriorate. This marks a limit of the
present surrogate. The experiment shows that boundary information is a
separate signal from logit distillation and exposes a Pareto trade-off
between accuracy, reference flips, boundary-band flips, and
boundary-mask preservation. Calibration boundary Jaccard remains
predictive of held-out boundary Jaccard for PTQ-W (\(r=0.947\)),
data-aware rounding-W (\(r=0.936\)), and boundary-aware rounding-W
(\(r=0.994\)).

Finally, we repeat the PTQ part of the protocol on full official
CIFAR-10 with a reduced residual CNN and three seeds. This full-dataset
scale check uses a compact residual architecture trained and evaluated
on the complete CIFAR-10 split. The mean full-precision held-out
accuracy is \(0.7000\) across seeds. Accuracy, reference flips, and
boundary-band flips are computed on the full held-out split; local
differential metrics use a stratified 3000-sample audit subset.

All entries are means over three seeds. The full-CIFAR run supports the
same qualitative conclusion at a larger scale. At 6 bits, PTQ-W loses
\(0.0029\) accuracy relative to float while changing \(5.3\%\) of
held-out decisions and \(24.5\%\) of low-margin boundary-band decisions.
PTQ-WA at 6 bits has similar accuracy loss and larger reference drift.
Accuracy stopping selects 6-bit quantization for both PTQ families in
all three seeds; stricter profiles select 8-bit or the full-precision
reference. Thus the boundary-aware criterion changes the engineering
decision even when accuracy appears acceptable.

\textbf{Reproducibility and data availability.} All datasets are public:
\texttt{sklearn\ digits}, MNIST, Fashion-MNIST, and the official
CIFAR-10 Python archive. The experiments are implemented as independent
scripts under \texttt{experiments/}. Raw metrics, stopping decisions,
bootstrap summaries, and calibration-to-test correlations are written as
CSV and JSON files in the results directory. The command lines for the
main CIFAR-10 PTQ/QAT run, the full-CIFAR reduced-residual scale check,
and the rounding ablations are provided in the supplementary
reproducibility file, together with the random seeds, bit-widths,
quantization families, audit-grid resolution, rounding-sample counts,
and CPU-thread settings used to generate the reported tables.

\FloatBarrier

\section{Conclusion}\label{conclusion}

We presented a boundary-aware diagnostic protocol for neural network
quantization. The protocol measures how a quantized classifier changes
the decision geometry of a full-precision reference through local
boundary displacement, normal variation, slice-boundary Jaccard
distance, grid prediction changes, and multiclass junction stability.
Experiments across several datasets and architectures show that low-bit
quantization and QAT can preserve or improve accuracy while changing the
reference boundary. Boundary-aware stopping uses calibration geometry to
select more conservative bit-widths and reduce held-out label flips.
Boundary-band analysis shows that quantization regressions concentrate
near low-margin reference samples, which is the region where deployment
updates are most likely to change decisions. An initial fixed-bit
rounding experiment shows that boundary information can be inserted into
the quantization objective, with a modest and tuning-dependent
trade-off. The central claim is that quantization can reconfigure
decision boundaries at deployment-relevant finite scales, and that this
behavior can be measured, predicted, and used as an explicit
model-selection criterion.

\section*{Declaration of Competing Interest}

The authors declare the absence of known competing financial interests
or personal relationships that could have appeared to influence the work
reported in this paper.

\FloatBarrier

\bibliographystyle{elsarticle-num}
\bibliography{references}

@inproceedings{szegedy2014intriguing,
  title = {Intriguing Properties of Neural Networks},
  author = {Szegedy, Christian and Zaremba, Wojciech and Sutskever, Ilya and Bruna, Joan and Erhan, Dumitru and Goodfellow, Ian and Fergus, Rob},
  booktitle = {International Conference on Learning Representations},
  year = {2014},
}

@inproceedings{goodfellow2015explaining,
  title = {Explaining and Harnessing Adversarial Examples},
  author = {Goodfellow, Ian J. and Shlens, Jonathon and Szegedy, Christian},
  booktitle = {International Conference on Learning Representations},
  year = {2015},
}

@inproceedings{moosavi2016deepfool,
  title = {DeepFool: A Simple and Accurate Method to Fool Deep Neural Networks},
  author = {Moosavi-Dezfooli, Seyed-Mohsen and Fawzi, Alhussein and Frossard, Pascal},
  booktitle = {Proceedings of the IEEE Conference on Computer Vision and Pattern Recognition},
  year = {2016},
  pages = {2574--2582},
}

@inproceedings{moosavi2019curvature,
  title = {Robustness via Curvature Regularization, and Vice Versa},
  author = {Moosavi-Dezfooli, Seyed-Mohsen and Fawzi, Alhussein and Uesato, Jonathan and Frossard, Pascal},
  booktitle = {Proceedings of the IEEE/CVF Conference on Computer Vision and Pattern Recognition},
  year = {2019},
  pages = {9078--9086},
}

@inproceedings{jacob2018quantization,
  title = {Quantization and Training of Neural Networks for Efficient Integer-Arithmetic-Only Inference},
  author = {Jacob, Benoit and Kligys, Skirmantas and Chen, Bo and Zhu, Menglong and Tang, Matthew and Howard, Andrew and Adam, Hartwig and Kalenichenko, Dmitry},
  booktitle = {Proceedings of the IEEE Conference on Computer Vision and Pattern Recognition},
  year = {2018},
  pages = {2704--2713},
}

@article{krishnamoorthi2018quantizing,
  title = {Quantizing Deep Convolutional Networks for Efficient Inference: A Whitepaper},
  author = {Krishnamoorthi, Raghuraman},
  journal = {arXiv preprint arXiv:1806.08342},
  year = {2018},
}

@inproceedings{nagel2020up,
  title = {Up or Down? Adaptive Rounding for Post-Training Quantization},
  author = {Nagel, Markus and Amjad, Rana Ali and van Baalen, Mart and Louizos, Christos and Blankevoort, Tijmen},
  booktitle = {Proceedings of the 37th International Conference on Machine Learning},
  year = {2020},
  pages = {7197--7206},
  series = {Proceedings of Machine Learning Research},
  volume = {119},
}

@article{mickisch2020understanding,
  title = {Understanding the Decision Boundary of Deep Neural Networks: An Empirical Study},
  author = {Mickisch, David and Assion, Felix and Gre{\ss}ner, Florens and G{\"u}nther, Wiebke and Motta, Mariele},
  journal = {arXiv preprint arXiv:2002.01810},
  year = {2020},
}

@article{lin2019defensive,
  title = {Defensive Quantization: When Efficiency Meets Robustness},
  author = {Lin, Ji and Gan, Chuang and Han, Song},
  journal = {arXiv preprint arXiv:1904.08444},
  year = {2019},
}

@article{gorsline2021adversarial,
  title = {On the Adversarial Robustness of Quantized Neural Networks},
  author = {Gorsline, Micah and Smith, James and Merkel, Cory},
  journal = {arXiv preprint arXiv:2105.00227},
  year = {2021},
}

@article{song2021layerwise,
  title = {A Layer-wise Adversarial-aware Quantization Optimization for Improving Robustness},
  author = {Song, Chang and Ranjan, Riya and Li, Hai},
  journal = {arXiv preprint arXiv:2110.12308},
  year = {2021},
}

@book{maclagan2015introduction,
  title = {Introduction to Tropical Geometry},
  author = {Maclagan, Diane and Sturmfels, Bernd},
  series = {Graduate Studies in Mathematics},
  volume = {161},
  publisher = {American Mathematical Society},
  address = {Providence, RI},
  year = {2015}
}

@incollection{richter2005first,
  title = {First Steps in Tropical Geometry},
  author = {Richter-Gebert, J{\"u}rgen and Sturmfels, Bernd and Theobald, Thorsten},
  booktitle = {Idempotent Mathematics and Mathematical Physics},
  editor = {Litvinov, Grigori L. and Maslov, Victor P.},
  series = {Contemporary Mathematics},
  volume = {377},
  pages = {289--317},
  publisher = {American Mathematical Society},
  address = {Providence, RI},
  year = {2005},
}

\end{document}